\documentclass{amsart}
\usepackage{url}
\usepackage{nohyperref}
\usepackage{graphics}
\begin{document}
\noindent
{\small
Topology Atlas Invited Contributions \textbf{7} (2002) no.~1, 3 pp.}
\bigskip
\title[Algebraic formulas imitating Pachner moves]{Algebraic formulas 
whose structure imitates Pachner moves and new types of acyclic complexes}
\author{Igor G. Korepanov}
\address{Southern Ural State University, 76 Lenin avenue, Chelyabinsk 
454080, Russia}
\email{kig@susu.ac.ru}
\maketitle

A triangulation of a manifold can be transformed into a different 
triangulation by a sequence of Pachner moves.\footnote{Let me proceed 
here in such easy style. Of course, I am omitting many technical details.}  
In the three-dimensional case, there are four of them: moves 
$2 \to 3$, $1 \to 4$ and two inverse moves. 
A move $2 \to 3$ replaces two adjacent tetrahedra $ABCD$ and $EABC$ with 
three tetrahedra $ABED$, $BCED$ and $CAED$ occupying the same domain and 
having the same common boundary.

\begin{figure}
\includegraphics{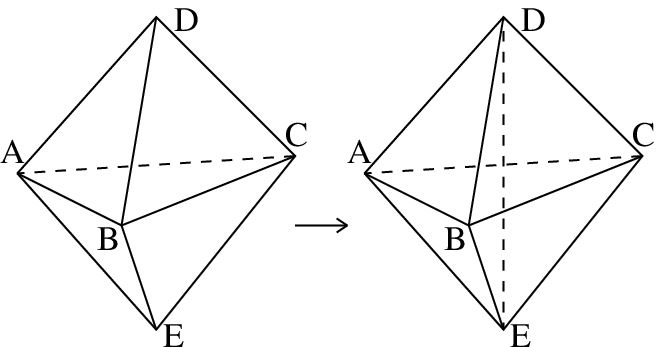}
\caption{A $2 \to 3$ Pachner move}
\end{figure}

A move $1 \to 4$ adds a new vertex $E$ inside a tetrahedron $ABCD$ and 
decomposes it into four tetrahedra $ABCE$, $ABED$, $AECD$ and $EBCD$.

\begin{figure}
\includegraphics{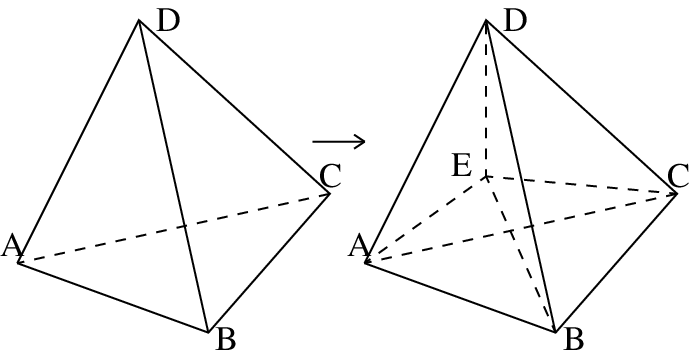}
\caption{A $1 \to 4$ Pachner move}
\end{figure}

Similar statements hold in higher dimensions as well. For instance, in 
four dimensions\footnote{I am of course speaking about the 
piecewise-linear category.} there are \emph{five} Pachner moves: 
$3 \to 3$, $2 \leftrightarrow 4$ and $1 \leftrightarrow  5$.

Suppose I want to invent some (new) manifold invariants. A natural way to 
do so is to invent an algebraic object which would correspond to a 
triangulation and whose changes under the Pachner moves of the 
triangulation would be describable in some simple way. This implies 
assigning some algebraic values to elements of the triangulation. For 
example, in papers \cite{K1, K2, K3, KM2} these are Euclidean 
metric  values: edge lengths, volumes of tetrahedra or $4$-simplices, 
dihedral angles, etc.; but a different example is presented in 
\cite{KM1}.

First, I need a \emph{local} formula whose structure would imitate a 
Pachner move. It turns out that a key role is played by moves $2 \to 3$ 
in three dimensions and by moves $3 \to 3$ in four dimensions. So, the 
local formula for three dimensions must contain values corresponding to 
$2$ initial tetrahedra in its l.h.s., and those corresponding to $3$ 
resulting tetrahedra---in its r.h.s. A good example of a formula of such 
kind is provided by the \emph{pentagon equation for quantum $6j$-symbols}. 
My formulas are, however, \emph{classical}. One of them reads:
$$V_{ABCD} V_{EABC} = - 6 V_{ABED} V_{BCED} V_{CAED} a / DE^2,$$
where $DE$ means the length of the edge $DE$ which is added to the complex 
when we pass from the two tetrahedra in the l.h.s.\ to three ones in the 
r.h.s., and $a$ is the partial derivative of the \emph{deficit angle}
around edge $DE$ (in the sense widely used in discrete gravity theories) 
with respect to its length.

This very formula can be obtained from the quantum pentagon equation using 
the quasiclassics of quantum $6j$-symbols,%
\footnote{Actually, it requires taking \emph{two} quasiclassics in the 
following different senses: first, the quantum parameter $q$ tends to its 
classical value $1$, which gives classical $6j$-symbols, and then the 
dimensions of the irreps go to infinity. About the latter quasiclassics 
see \cite{R} and references therein to works by Regge and Ponzano. These 
quasiclassics enable one to obtain the above formula using the stationary 
phase method.} 
but nothing like that is known for its four-dimensional analogue from 
\cite{K2} and its $SL(2)$-analogue from \cite{KM1}. This poses at once 
the \emph{first} (and most interesting for me) \emph{question}: 
what are the quantum analogues of other classical formulas  (which still 
have a \emph{very quasiclassical appearance})?

Second, there must exist \emph{global} algebraic objects corresponding to 
the whole manifold (not just to a cluster of neighbouring simplices). 
Here I arrived at the Jacobian matrix $A$ of partial derivatives of 
\emph{all} deficit angles in \emph{all} edge lengths (and its analogues 
for other cases). 
Then, after some two years of thinking, I realized that $A$ must be 
regarded as one of the linear mappings in an acyclic complex of based 
vector spaces! 
That kind of acyclic complex looks pretty unusual: its vector spaces 
consist, in particular, of infinitesimal deformations of algebraic (e.g., 
Euclidean metric) values associated with the simplicial complex. See my 
recent paper [K3] for the description of acyclic complex for some 
three-dimensional manifolds and its partial description for 
four-dimensional manifolds. The invariant of \emph{all} Pachner moves is 
expressed through the \emph{torsion} of that complex, and this poses the 
\emph{second question}: 
how can we include the \emph{Reidemeister} torsion in this scheme?

So, the expected result of this work must be some theory uniting together 
such subjects as different analogues of $6j$-symbols, torsion of acyclic 
complexes and hopefully other interesting things. I invite my colleagues 
throughout the world to work together in this direction.


\begin{thebibliography}{9}

\bibitem{K1}
Korepanov, I.~G., Invariants of PL manifolds from metrized simplicial 
complexes. Three-dimensional case, J.\ Nonlin.\ Math.\ Phys. {\bf 8}:2 
(2001), 196--210. 
arXiv:math.GT/0009225.
\url{http://www.sm.luth.se/~norbert/home_journal/electronic/82lett2.pdf}

\bibitem{K2}
Korepanov, I.~G., Euclidean 4-simplices and invariants of four-dimensional 
manifolds. I. Moves $3 \to 3$,  Accepted for publication in 
Theor.\ Math.\ Phys. 
arXiv:math.GT/0211165.
\url{http://www.lib.csu.ru/ik/3-3eng.pdf}

\bibitem{K3}
Korepanov, I.~G., Euclidean $4$-simplices and invariants of 
four-dimensional manifolds. II. An algebraic complex and moves 
$2 \leftrightarrow 4$, Submitted to Theor.\ Math.\ Phys. 
arXiv:math.GT/0211166
\url{http://www.lib.csu.ru/ik/acycleng.pdf}

\bibitem{KM1}
Korepanov, I.~G., and Martyushev, E.~V., A classical solution of the 
pentagon equation related to the group $SL(2)$, Theor.\ Math.\ Phys. 
{\bf 129}:1 (2001), 1320-1324. 

\bibitem{KM2}
Korepanov, I.~G., and Martyushev, E.~V., Distinguishing three-dimensional 
lens spaces $L(7,1)$ and $L(7,2)$ by means of classical pentagon equation, 
J.\ Nonlin.\ Math.\ Phys. {\bf 9}:1 (2002), 86--98. 
arXiv:math.GT/0210343.
\url{http://www.sm.luth.se/~norbert/home_journal/electronic/91art5.pdf}

\bibitem{R}
Roberts, J.~D., Classical $6j$-symbols and the tetrahedron, Geometry and 
Topology {\bf 3} (1999), 21--66. 
arXiv:math-ph/9812013.
\url{http://www.maths.warwick.ac.uk/gt/GTVol3/paper2.abs.html}

\end{thebibliography}
\end{document}